\documentclass[a4paper,12pt]{article}
\usepackage{amsmath, amsfonts, amssymb}
\usepackage[english, russian]{babel}
\begin{document}

\textit{}
\begin{center}\textbf{On absolute nilpotent and idempotent elements of an evolution algebra corresponding to permutations.}

\medskip

\textbf{\verb"B. A. Narkuziev"}\\
\textit{Samarqand State University, Samarqand, Uzbekistan}\\
\smallskip
\textit{e-mail: bnarkuziev@yandex.ru}
\end{center}

\begin{quote}
\textbf{Abstract.}We describe  absolute nilpotent and some idempotent elements of an $n$- dimensional evolution algebra corresponding to two permutations and we decompose such algebras to the direct sum of evolution algebras corresponding to cycles of the permutations.

\noindent{\bf Keywords:} Evolution algebra, algebra of permutations,absolute nilpotent and idempotent element.
{\small }
\end{quote}

\medskip
\textbf{1. Introduction}

\

An evolution algebra is an abstract system, it gives an insight for the study of non-Mendelian genetics. Results on genetic evolution and genetic algebras, can be found in Lyubich's book \cite{Lyu1992}.
In Tian's book \cite{Tian2008} the foundations of evolution algebras are developed and several basic properties are studied. The concept of evolution algebra lies between algebras and dynamical systems. Algebraically, evolution algebras are non-associative Banach algebras; dynamically, they represent discrete dynamical systems. Evolution algebras have many connections with other mathematical fields including graph theory, group theory, Markov chains, dynamic systems, stochastic processes, mathematical physics, etc.\cite{CasLad2011}-\cite{Lyu1992}.

In this paper we describe absolute nilpotent and idempotent elements of an evolution algebra constructed by two permutations.

The paper is organized as follows. Section 2 is devoted to absolute nilpotent elements of evolution algebras corresponding to permutations and to a criteria for this algebra to be baric. In Section 3 we study idempotent elements of the  evolution algebra  and give a general  analysis of idempotent elements of the two- dimensional evolution algebra corresponding to permutations. Section 4 is devoted to the isomorphism of evolution algebra corresponding to  permutations.

\

\textbf{2. Absolutely nilpotent elements.}

\

\textbf{Definition 2.1.} Let $(E,\cdot)$ be an  algebra  over  a  field  $K$. If it admits a basis $e_1,e_2...,$ such that
 $$e_{i}\cdot e_{j}=0, \ \ \ \ \ \ \ \ \ \ if \ \ \ i\neq j \ ; $$
$$
\ e_{i}\cdot e_{i}=\sum\limits_{k} a_{ik}e_k \ ,\ \ \ \ \ \ \ \ \ \
 for\ \ \ any\ \ \ \ i\ ,
$$
 then this algebra is called an $evolution \  algebra$. This basis is  called natural basis.

We denote by $M=(a_{ij})$  the  matrix  of    the structural constants of the evolution algebra  $E$.

Let $S_{n}$ be  the group  of   permutations  of  degree   $n$. Take $\pi,\tau \in S_n,$

$$ \pi=\left( \begin{array}{cc} 1 \,\,\,\ \ \ \ 2 \ \ \ \,\,\,\ 3 \ \ \ \ \ \ \ ...\ \ \ \,\,\,\ n \\
\pi(1) \,\ \pi(2) \,\ \pi(3) \ \ \  ... \,\ \ \ \pi(n)\end{array} \right),\ \ \ $$

$$ \ \ \ \ \ \ \ \ \ \ \ \ \ \ \ \ \ \ \ \tau=\left( \begin{array}{cc} 1 \,\,\,\ \ \ \ 2 \ \ \ \,\,\,\ 3 \ \ \ \ \ \ \ ...\ \ \ \,\,\,\ n \\
\tau(1) \,\ \tau(2) \,\ \tau(3) \ \ \  ... \,\ \ \ \tau(n)\end{array} \right), \ \ \ \pi(i), \tau(i)\in \{1,...,n\}.$$

Consider a $n$-dimensional evolution  algebra  $E^{n}_{\pi,\tau}$ over  the  field   $R$  with  a  finite  natural  basis   $e_1,e_2,...e_n,$ and
    multiplication   given  by:
    $$e_{i}\cdot e_{j}=0, \ \  if \ \  i\neq j \ ;\ \ \ \ $$
 $$  e_{i}\cdot e_{i}=a_{i\pi(i)}\cdot e_{\pi(i)} + a_{i\tau(i)}\cdot e_{\tau(i)},\ \ \
 for\ \ \ any\ \ \ \ i\ ,$$
then this algebra  is  called  an  \emph{evolution algebra corresponding to permutations} $\pi$ and $\tau$.

We  note  that $a_{ij}=\left\{\begin{array}{llllll} a_{i\pi(i)}, \ \ \ \ \  \mbox{if} \ \ \ j=\pi(i)\ \\
a_{i\tau(i)}, \ \ \ \ \  \mbox{if} \ \ \ j=\tau(i)\ \\
0\ \ \  ,\,\ \ \ \ \ \ \ \mbox{if}\ \ \  j\neq \pi(i),\tau(i) \end{array} \right. \ \ \ \ \ \ i={1,...,n}  .$

We assume that $ \pi\neq\tau.$ For the case $ \pi=\tau$ some properties of the algebra are known \cite{Nar2014}-\cite{KhudOmirov2015}.

The following properties  of  evolution algebras are known \cite{Tian2008}:

(1) Evolution algebras are not associative, in general.

(2) Evolution algebras are commutative, flexible.

(3) Evolution algebras are not power-associative, in general.

(4) The direct sum of evolution algebras is also an evolution algebra.

(5) The Kronecker product of evolution algebras is an evolution algebra.

These  properties  are  hold  for  $E^{n}_{\pi,\tau}$   too.

A  $character$  for  an  algebra   $\mathcal A$  is  a  nonzero  multiplicative  linear form on  $\mathcal A$ , that is,
 a   nonzero  algebra  homomorphism  from   $\mathcal A$  to R \cite{Lyu1992} .\\

\textbf{Definition 2.2.}  A   pair  $( \mathcal A ,\sigma)$  consisting \  of \  an\   algebra \  $\mathcal A$ \  and \  a \  character $\sigma$
 on \  $\mathcal A$  is  called   a\   baric\  algebra.\  The \ homomorphism $\sigma$  is\  called\   the\   weight (or baric) function
  of  $\mathcal A$  and  $\sigma(x)$ the weight (baric value) of  $x$.
The  following  theorem  is  a  particular  case  of  Theorem 3.2  in  \cite{CasLad2011}.

\

\textbf{Theorem 2.3.}  \emph{An  evolution algebra $ E^{n}_{\pi,\tau }$  over the field  R, is baric if and only if one of the following conditions holds:}

\emph{ 1) $k_{0}$ is a fixed point for $\pi$ (or $\tau $),$ \ \tau({k_{0}})\neq k_{0} \ (\pi(k_{0})\neq k_{0})$  and  its  matrix  $M=(a_{ij})\ \ i,j=1,...,n$ of structural constants satisfies  $a_{k_{0}k_{0}}\neq0$, $a_{\tau^{-1}(k_{0})k_{0}}=0$ (or $a_{\pi^{-1}(k_{0})k_{0}}=0$).
  In this case the corresponding weight function is $ \sigma(x)=a_{k_{0}k_{0}}\cdot x_{k_{0}}$ .}

  \emph{ 2) $k_{0}$ is a fixed point for both $\pi$ and $\tau$, and $a_{k_{0}k_{0}}\neq0.$ In this case the corresponding weight function is $ \sigma(x)=2\cdot a_{k_{0}k_{0}}\cdot x_{k_{0}}$ .}

\medskip
\textbf{Corollary 2.4.}  \emph{If the permutations $\pi$ and $ \tau$, mentioned in Theorem 2.3 have $m$ fixed points $k_{m}, \ m\leq n$ which satisfy conditions of Theorem 2.3, then the evolution algebra $ E^{n}_{\pi,\tau }$ has exactly m weight functions $ \sigma(x)=a_{k_{m}k_{m}}\cdot x_{k_{m}} $ (or $ \sigma(x)=2a_{k_{m}k_{m}}\cdot x_{k_{m}}$ ).}\\

  \textbf{Definition 2.5.} An element $x\in E$ is called an absolute nilpotent if $x^{2}=0$.

  Let $E=R^{n}$ be an evolution algebra over the  field $R$ with the matrix
   \\ $M=(a_{ij})$ of structural constants, then for arbitrary $x=\sum\limits_{i}x_{i}e_{i}$ and
   $y=\sum\limits_{i}y_{i}e_{i}\in R^{n}$ we have $x\cdot y=\sum\limits_{i}x_{i}y_{i}e_{i}e_{i}$ , $ x^{2}=\sum\limits_{i}x^{2}_{i}e_{i}e_{i}$ and for $ E^{n}_{\pi,\tau }$ we have
   $$ x^{2}=\sum \limits_{i}x^{2}_{i}(a_{i\pi(i)}\cdot e_{\pi(i)} + a_{i\tau(i)}\cdot e_{\tau(i)}).\ \ \ \ \ \ \ \ \  \ \ \ (2.1)$$
   Let $j_{k}=\tau^{-1}(\pi(k)),$ then $\sum \limits_{i}a_{i\tau(i)} x^{2}_{i}  e_{\tau(i)}=\sum \limits_{k} a_{j_{k}\pi(k)}x^2_{j_{k}}e_{\pi(k)} $  therefore from (2.1) we have
   $$ x^{2}=\sum \limits_{k=1}^{n}e_{\pi(k)}(a_{k\pi(k)}x^2_{k}+ a_{j_{k}\pi(k)}x^2_{j_{k}} ). $$
   Thus the equation $x^{2}=0$ will be $$ a_{k\pi(k)}x^2_{k}+ a_{j_{k}\pi(k)}x^2_{j_{k}}=0 ,\ \  k\in \{1,...,n\}.\ \ \ \ \ \ \ \ \ \ \ (2.2) $$
  If $detM\neq 0$ then the system (2.2) has unique solution $(0,...,0)$.
   If $detM=0$ and $rank$($M$)=$n-1$ then we can assume that the first $n-1$ rows of $M$ are linearly independent and denote
   $M_{n-1}=(a_{s,t})_{s,t=1,...,n-1}$\ \ ,
  $$M_{in}=\left( \begin{array}{cccc} a_{11}\ \ \ \  ...\ \ \ \  a_{i-1,1}\ \ \ \ \ \ a_{n1} \ \ \ \ \  a_{i+1,1} \ \ ... \ \ \  a_{n-1,1}
    \\ a_{12}\ \ \ \ ... \ \ \ \  a_{i-1,2}\ \ \ \ \ \ a_{n2} \ \ \ \ \ a_{i+1,2} \ \ ... \ \ \  a_{n-1,2} \\ ...\ \ \ \ \ \ \ \ \ \ \ \  ...\ \ \ \ \ \ \  \ \ \ \  ... \ \ \ \ \ \ \ \ \  ....\ \ \ \ \ \ \ \ \ \ ... \\
 a_{1,n-1}\ ... \ \ a_{i-1,n-1}\ \ a_{n,n-1} \ \ a_{i+1,r} \ ...  \ \ a_{n-1,n-1} \end{array} \right) \  .$$ \\

 \textbf{Proposition 2.6.} \  \emph{Let $M$ be the matrix of structural constants for $ E^{n}_{\pi,\tau }$.
  The finite dimensional evolution algebra $ E^{n}_{\pi,\tau }$ has the unique absolute nilpotent element $(0,...,0)$
  if one of the following conditions is satisfied.
 \\ (i) $det(M)\neq 0$ ;
 \\ (ii) $det(M)=0$, $rank(M)=n-1 $ and \ \  $det(M_{i_{0}n})\cdot  det(M_{n-1})>0$ \ \ for some $i_{0}\in \{1,...,n-1\}$;
 \\ (iii) For all $k\in \{1,...,n \}, \ a_{k\pi(k)}\cdot a_{j_{k}\pi(k)}>0 $ .}\\

\textbf{Proof.} \ The first and second parts of the proposition are particular cases of Proposition 4.12 in \cite {CasLad2011}.
Part 3 is simple consequence of (2.2).\ \

 To see difference between conditions $(ii)$ and $(iii)$ we give the  following examples.\\

 \textbf{Example 1.} \ (For the case $(ii)$ ) Let
 $$ \pi=\left( \begin{array}{cc} 1 \ \ \ 2 \ \ \ 3 \ \ \ 4 \\
 3 \ \ \ 1 \ \ \ 4 \ \ \ 2\end{array} \right),\ \ \tau=\left( \begin{array}{cc} 1 \ \ \ 2 \ \ \ 3 \ \ \ 4 \\
 2 \ \ \ 3 \ \ \ 4 \ \ \ 1\end{array} \right),\ \ \ $$ from the system (2.2) we have the following  system of equations
 $$ \left\{\begin{array}{llll} a_{13}x^2_{1}+ a_{23}x^2_{2}=0 \\
a_{21}x^2_{2}+ a_{41}x^2_{4}=0 \\
a_{34}x^2_{3}+ a_{34}x^2_{3}=0  \\
a_{42}x^2_{4}+ a_{12}x^2_{1}=0 \
\end{array} \right. \ \ \ and \ \ \ M = \left( \begin{array}{cccc} a_{13} \ \ \ \  a_{23} \ \ \  0 \ \ \ \ \ 0\\
  \ \  0 \ \ \ \ \  a_{21} \ \ \  \ 0 \  \ \ \ a_{41} \\
 0 \ \ \ \ \ \  0 \ \ \ \   2a_{34} \ \  0 \\
  \ \  a_{12} \ \ \ \ 0 \ \ \ \ \ \ 0 \ \ \ \ a_{42}\end{array} \right) \   $$
$ det M=2(a_{13}\cdot a_{21}\cdot a_{34}\cdot a_{42}+a_{12}\cdot a_{23}\cdot a_{34}\cdot a_{41}).$
If $a_{13}=-1$ and $a_{ij}=1,\ ij\neq 13 $ then
 $ M = \left( \begin{array}{cccc}  -1  \ \ \ \  1  \ \ \ \  0 \ \ \ \ 0 \\
 \ 0 \ \ \ \ \ 1  \ \ \  \ 0  \ \ \ 1 \\
 \ 0 \ \ \ \ \ 0 \ \ \ \   2  \ \ \ 0 \\
 \  1 \ \ \ \ \  0 \ \ \ \   0 \ \ \  1 \end{array} \right) $ and  $ det M=0$, $rank M=3 $.
 $M_{n-1}=M_{3}= \left( \begin{array}{ccc}  -1  \ \ \ \  1  \ \ \ \  0  \\
 \ 0 \ \ \ \ \ 1  \ \ \  \ 0   \\
 \ 0 \ \ \ \ \ 0 \ \ \ \   2    \end{array} \right),$
  $M_{i_{0}n}=M_{14}=  \left( \begin{array}{ccc}  0  \ \ \ \  1  \ \ \ \  0  \\
 1 \ \ \ \  1  \ \ \  \ 0   \\
 0 \ \ \ \ 0 \ \ \ \   2    \end{array} \right)$ ,\ \ \ \ \ \ \ \ \ \ \ \ \ \ \ \ \ \ \ \ \
 $ det M_{3} \cdot det M_{14}=(-2) \cdot (-2)=4$ .
 It is easy to see that the system $ \left\{\begin{array}{llll} -x^2_{1}+ x^2_{2}=0 \\
\ \ x^2_{2}+ x^2_{4}=0 \\
\ \ x^2_{3}+ x^2_{3}=0  \\
\ \ x^2_{4}+ x^2_{1}=0 \
\end{array} \right.$ has only trivial solution.

  \textbf{Example 2.} \ (For the case $(iii)$) Let
 $$ \pi=\left( \begin{array}{cc} 1 \ \ \ 2 \ \ \ 3 \ \ \ 4 \\
 3 \ \ \ 2 \ \ \ 4 \ \ \ 1\end{array} \right),\ \ \tau=\left( \begin{array}{cc} 1 \ \ \ 2 \ \ \ 3 \ \ \ 4 \\
 2 \ \ \ 3 \ \ \ 1 \ \ \ 4\end{array} \right),\ \ \ $$ from the system (2.2) we have the following  system of equations
 $$ \left\{\begin{array}{llll} a_{13}x^2_{1}+ a_{23}x^2_{2}=0 \\
a_{22}x^2_{2}+ a_{12}x^2_{1}=0 \\
a_{34}x^2_{3}+ a_{44}x^2_{4}=0  \\
a_{41}x^2_{4}+ a_{31}x^2_{3}=0 \
\end{array} \right.  \ \ \ and \ \ \  M = \left( \begin{array}{cccc} a_{13} \ \ \ \  a_{23} \ \ \ \  0 \ \ \ \ \ 0\\
  a_{12} \ \ \ \ a_{22} \ \ \  \ 0 \  \ \ \ \ 0  \\
  \ \ \ 0 \ \ \ \ \ \ \  0 \ \ \ \ \ a_{34} \ \ \ \ a_{44}  \\
  \ \ \ 0 \ \ \ \ \ \ \  0 \ \ \ \ \ a_{31} \ \ \ \ a_{41} \end{array} \right) \   $$
if all $a_{ij}=1$ then $det M =0,rank M =2$ and this system has only trivial solution.

  Let $$ \tau^{-1}=\left( \begin{array}{cc} 1 \,\,\,\ \ \ \ \ \  \ 2 \ \ \  \  \ \,\, ...\ \ \ \ \ \,\,\ n \\
\tau^{-1}(1) \,  \ \tau^{-1}(2) \, \ \  ... \,\  \tau^{-1}(n)\end{array} \right),\  $$ and $$ \pi \tau^{-1} =\left( \begin{array}{cc} 1 \,\,\,\ \ \ \ 2  \ \ \ \ \ \ ...\ \ \ \,\,\,\ n \\
j_{1} \,\ \ \ \  j_{2} \,\  \ \ \  ... \,\ \ \ \ \ \ j_{n}\end{array} \right),\ \ \ \ \ where \ \ j_{k}=\tau^{-1}(\pi(k)).$$

The following theorem and proposition fully describes absolute nilpotent elements of $ E^{n}_{\pi,\tau }$.

\

  \textbf{Theorem 2.7.}  \emph{ Let $x^{*}=(x^{*}_{1},x^{*}_{2},...,x^{*}_{n})$ be an absolute nilpotent element for $ E^{n}_{\pi,\tau }$ and let
  $(l_{1},l_{2},...,l_{p})$ be one of the independent cycles of $\pi \tau^{-1},\ (p\leq n).$\\
    1) If $a_{l_{k}\pi(l_{k})}\cdot a_{l_{k}\tau (l_{k})}\neq 0, $ for all $k=1,2,...,p,$ then $x^{*}_{l_{1}}=x^{*}_{l_{2}}=...=x^{*}_{l_{p}}=0$ or
$\prod \limits_{i=1}^{p}x^{*}_{l_{i}}\neq 0.$\\
 2) If $a_{l_{s}\pi(l_{s})}=0$ , $a_{l_{s}\tau(l_{s})}=0$ for some $s$ and $a_{l_{k}\pi(l_{k})}\cdot a_{l_{k}\tau (l_{k})}\neq 0,$ for any
$k\neq s, \ k=1,...,p$ then $x^{*}_{l_{s}}$ will be any real number and $x^{*}_{l_{i}}=0, \ i\neq s , \ i=1,...,p.$ }\\
 3) If  $a_{l_{s}\pi(l_{s})}=0$ (or $a_{l_{s}\tau(l_{s})}=0$ ) for some $s$  $( 1 \leq s \leq p)$  and $a_{l_{m}\tau (l_{m})}\neq 0,$ for any $m=1,...,p$ ( or $a_{l_{m}\pi (l_{m})}\neq 0,$ for any $m=1,...,p$ ) then $x^{*}_{l_{1}}=x^{*}_{l_{2}}=...=x^{*}_{l_{p}}=0$ .

\

  \textbf{Proof.}  \ 1)\ From the system of equations (2.2) we have
$$ a_{l_{1}\pi(l_{1})}x^2_{l_{1}}+ a_{j_{l_{1}}\pi(l_{1})}x^2_{j_{l_{1}}}=0 .\ \ \ \ \ \ \ \ \ \ \ \ \ \ \ \ \ \  (2.3)  $$
 We know $j_{k}=\tau^{-1}(\pi(k)) $ it means that $l_{k}$ and $j_{l_{k}}$ are in one cycle of $\pi\tau^{-1}$, also
   $$j_{l_{k}}=l_{k+1},\  \pi(l_{k})=\tau(l_{k+1})\ \ \ \ \ \ \ \ \ \ \ \ \ \ \ \ \ \ \ \ (2.4)$$
   similarly  $j_{l_{1}}=\tau^{-1}(\pi(l_{1})),\
 \tau(j_{l_{1}})=\pi(l_{1}), \ \tau(l_{2})=\pi(l_{1})$ thus equation (2.3) will be $$ a_{l_{1}\pi(l_{1})}x^2_{l_{1}} + a_{l_{2}\tau(l_{2})}x^2_{l_{2}}=0. $$
  Since $a_{l_{1}\pi(l_{1})}\neq 0,\ a_{l_{2}\tau(l_{2})}\neq 0$, if  $x_{l_{1}}=0,$ then $x_{l_{2}}=0$ too. By using method of mathematical induction it is easy to see that if
 $x_{l_{1}}=0 $ then $x_{l_{2}}=x_{l_{3}}=...=x_{l_{p}}=0.$

  2) From the system of equations (2.2) we obtain
  $$ \left\{\begin{array}{llll} a_{l_{1}\pi(l_{1})}x^2_{l_{1}}+ a_{j_{l_{1}}\pi(l_{1})}x^2_{j_{l_{1}}}=0 \\
a_{l_{2}\pi(l_{2})}x^2_{l_{2}}+ a_{j_{l_{2}}\pi(l_{2})}x^2_{j_{l_{2}}}=0 \\
....................................... \ \ \ \ \ \ \ \ \ \ \ \  \ \ \ \ \ \ \ \ \ \ \ \  (2.5)\\
a_{l_{p}\pi(l_{p})}x^2_{l_{p}}+ a_{j_{l_{p}} \pi(l_{p})}x^2_{j_{l_{p}}}=0
\end{array} \right. $$
by using equality (2.4) and $l_{p+1}=l_{1}$ the system (2.5) will be
$$ \left\{\begin{array}{llll} a_{l_{1}\pi(l_{1})}x^2_{l_{1}}+ a_{l_{2}\tau(l_{2})}x^2_{l_{2}}=0 \\
a_{l_{2}\pi(l_{2})}x^2_{l_{2}}+ a_{l_{3}\tau(l_{3})}x^2_{l_{3}}=0 \\
....................................... \\
a_{l_{p}\pi(l_{p})}x^2_{l_{p}}+ a_{l_{1}\tau(l_{1})}x^2_{l_{1}}=0 \ \ \ \ \ \ \ \ \ \ \ \ \ \ \ \ \ \ \ \ \ \ \ \  (2.6)
 \end{array} \right. $$
without loss of generality assume that  $a_{l_{1}\pi(l_{1})}=0, \ a_{l_{1}\tau(l_{1})}=0$ and $a_{l_{k}\pi(l_{k})}\cdot a_{l_{k}\tau(l_{k})}\neq 0,\ \ k=2,...,p $ then from the system of equation (2.6) it is easy to see that $x_{l_{2}}=x_{l_{3}}=...=x_{l_{p}}=0$ and $x_{l_{1}}$ is an arbitrary real number.\\
  3) Proof of this case is simple consequence of (2.6)

\

\textbf{Proposition 2.8.} \
1) \emph{If $a_{l_{k}\pi(l_{k})}\cdot a_{l_{k+1}\tau (l_{k+1})}\neq 0, $ for all $k=1,2,...,p,$ and $a_{l_{k_{0}}\pi(l_{k_{0}})} \cdot a_{l_{k_{0}+1}\tau(l_{k_{0}+1})} >0 $ for some $k_{0}=1,...,p$
then $x^{*}_{l_{1}}=x^{*}_{l_{2}}=...=x^{*}_{l_{p}}=0$, where $l_{p+1}=l_{1}$. \\
2) If $a_{l_{k}\pi(l_{k})} \cdot a_{l_{k+1}\tau(l_{k+1})} <0 $ for all $k=1,...,p $ then
$$ |x^{*}_{l_{k}}|=\sqrt{\frac{(-1)^{k-1}a_{l_{1}\pi(l_{1})}a_{l_{2}\pi(l_{2})}\cdot ... \cdot a_{l_{k-1}\pi(l_{k-1})}}{a_{l_{2}\tau (l_{2})}a_{l_{3}\tau (l_{3})}\cdot ... \cdot a_{l_{k}\tau (l_{k})}}} \cdot |x^{*}_{l_{1}}| . \  $$ }

\textbf{Proof.} 1)\ Using the equality (2.4) the system of equations (2.2) may be write as follows
$$ a_{l_{k}\pi(l_{k})}x^2_{l_{k}}+ a_{l_{k+1}\tau(l_{k+1})}x^2_{l_{k+1}}=0 . $$
If $a_{l_{k_{0}}\pi(l_{k_{0}})}\cdot a_{l_{k_{0}+1}\tau(l_{k_{0}+1})} >0$ for some $ k_{o}=1,...,p $ then
$ x^{*}_{l_{k_{0}}}=x^{*}_{l_{k_{0}+1}}=0$, by Theorem 2.7  we have $x^{*}_{l_{1}}=x^{*}_{l_{2}}=...=x^{*}_{l_{p}}=0$. \\
2)  If  $a_{l_{k}\pi(l_{k})} \cdot a_{l_{k+1}\tau(l_{k+1})} <0$ for all $k=1,...,p $ then nontrivial solutions of the system (2.6) are
$$ |x_{l_{2}}|=  \sqrt{-\frac{a_{l_{1}\pi(l_{1})}}{a_{l_{2}\tau(l_{2})}}} \cdot |x_{l_{1}}| ,$$
$$ |x_{l_{3}}|=  \sqrt{\frac{a_{l_{1}\pi(l_{1})} \cdot a_{l_{2}\pi(l_{2})}}{a_{l_{2}\tau(l_{2})}a_{l_{3}\tau(l_{3})}}} \cdot |x_{l_{1}}|, $$ $$  ..............................................$$
$$ |x_{l_{p}}|=\sqrt{\frac{(-1)^{p-1}a_{l_{1}\pi(l_{1})}a_{l_{2}\pi(l_{2})}\cdot ... \cdot a_{l_{p-1}\pi(l_{p-1})}}{a_{l_{2}\tau (l_{2})}a_{l_{3}\tau (l_{3})}\cdot ... \cdot a_{l_{p}\tau (l_{p})}}} \cdot |x_{l_{1}}|. $$

For an arbitrary $n$-dimensional  evolution algebra $ E $ over the field $R$ with the matrix $M=(b_{ij})$ of structural constants, the equation $x^{2}=0$
is given by the following system $$\sum\limits_{i}b_{ij}x^{2}_{i}=0,\ \ \ \ \ \ i,j=1,...,n \ \ \ \ \ \ \ \ \ \ \ \ \ \ \ \ \ \ (2.7) $$
If $det(M)\neq 0$ then the system $(2.7)$ has the unique solution $(0,...,0).$ If $det(M)= 0$ and $rank(M)=r$ then we can assume that the first $r$ rows of $M$
are linearly independent, consequently, the system $(2.7)$ can be written as
$$ x^{2}_{i}=-\sum \limits_{j=r+1}^{n}d_{ij}\cdot x^{2}_{j},\ \ \ \ \ \ \ \ i=1,...,r \ \ \ \ \ \ \ \ \ \ \ \ \ \ \ \ \ \  (2.8) $$
where $$d_{ij}=\frac{det(M_{ij})}{det(M_{r})},\  M_{r}=(b_{i_{*},j_{*}})_{i_{*},j_{*}=1,...,r}, \ \ $$
$$ M_{ij}= \left( \begin{array}{cccc} b_{11}\ ... \ b_{i-1,1}\ b_{j1} \ b_{i+1,1} \ ...  \ b_{r1}
    \\ b_{12}\ ... \ b_{i-1,2}\ b_{j2} \ b_{i+1,2} \ ...  \ b_{r2} \\ ...\ \ \ \ \ \ ... \ \ \ \ \ \ ... \ \ \ \ \ \ \ .... \\
 b_{1r}\ ... \ b_{i-1,r}\ b_{jr} \ b_{i+1,r} \ ...  \ b_{rr} \end{array} \right) \  . $$

An interesting problem is to find a necessary and sufficient condition on matrix $ D=(d_{ij})_{i=1,...,r ;  \ j=r+1,...,n }\ \ $ under which the system (2.8) has a unique solution. Here we shall
consider the case $rank(M)=n-2 $ (In \cite{CasLad2011} the case $rank(M)=n-1 $ is considered ).
\\ Then the system eq.(2.8) will be $$ x^{2}_{i}=-(d_{i,n-1}\cdot x^{2}_{n-1}+ d_{i,n}\cdot x^{2}_{n}),\ \ i=1,...,n-2.\ \ \ \ \ \ \ \ \  (2.9) $$

 \textbf{Proposition 2.9.}\ \emph{If $det(M)=0$ and $rank(M)=n-2$ then the evolution algebra $E=R^{n}$ has a unique absolutely nilpotent element $(0,...,0)$ if and only if
  $$ \left\{\begin{array}{ll} d_{i_{0},n-1}>0 \\ d_{i_{0},n}>0
\end{array} \right.  \ \ \ \ \ \ \ \ \ \ \ \ \ \ \ \ \ \ \ \ \ \ \ \ \ \ \ \ \ \ \ \ \\ \ \ \ \ \ \ \ \ \ \ \ \ \ \ \   (2.10)$$
for \ some\  $ i_{0}=1,...,n-2   $.}\\

\textbf{Proof}(\textit{Necessity}). Let  $ d_{i_{0},n-1}>0, d_{i_{0},n}>0$ for some $i_{0}$ , then solutions of  equation $ x^{2}_{i_{0}}+ d_{i_{0},n-1}\cdot x^{2}_{n-1}+ d_{i_{0},n}\cdot x^{2}_{n}=0$ are $ x_{i_{0}}=x_{n-1}=x_{n}=0$.It means that  the system of equations (2.9) has a unique trivial solution.\\
(\textit{Sufficiency}).Assume that the condition (2.10) is not satisfied for any $i_{0}$, then it is enough to consider the following cases
\\(A) Let $ \left\{\begin{array}{ll} d_{i,n-1}\leq 0 \\ d_{i,n}\leq 0
\end{array} \right. $ \ \  for \ any\ $ i=1,...,n-2  $ then $-(d_{i,n-1}\cdot x^{2}_{n-1}+ d_{i,n}\cdot x^{2}_{n})\geq 0 $ and we can define $x_{i}=\sqrt{-(d_{i,n-1}\cdot x^{2}_{n-1}+ d_{i,n}\cdot x^{2}_{n})}, \ i=1,...,n-2 $ in that case solutions of the system eq.(2.9) are
$$ (\sqrt{-(d_{1,n-1}\cdot x^{2}_{n-1}+ d_{1,n}\cdot x^{2}_{n})},...,\sqrt{-(d_{n-2,n-1}\cdot x^{2}_{n-1}+ d_{n-2,n}\cdot x^{2}_{n})},x_{n-1},x_{n})$$
where $ x_{n-1},x_{n}$ are any real numbers.\\
(B) Without loss of generality assume that $ \left\{\begin{array}{ll} d_{i,n-1}\geq 0 \\ d_{i,n}\leq 0
\end{array} \right. $ \ for \ any\ $ i=1,...,n-2 . $ In this case one of the nontrivial solutions of the system eq.(2.9) is
$$ x_{i}=0, i=1,...,n-2,\  x_{n-1}=\sqrt{-d_{i,n}}, \ x_{n}=\sqrt{d_{i,n-1}}.$$ \\

\textbf{3. Idempotent elements.}

\

 \textbf{Definition 3.1.} An element $x\in E$ is called $idempotent $ if $x^{2}=x$.
 \\ Such elements of an evolution algebra are especially important, because they are the fixed points of the evolution map $V(x)=x^{2}$, i.e. $V(x)=x$.

 For  $ E^{n}_{\pi,\tau }$ we have $x=\sum\limits_{i}x_{i}e_{i}= \sum\limits_{k}x_{\pi(k)}e_{\pi(k)}$, and the equation $x^{2}=x$ will be (see eq(2.2))
 $$ a_{k\pi(k)}x^2_{k}+ a_{j_{k}\pi(k)}x^2_{j_{k}}=x_{\pi (k)}, \ \  k \in \{1,...,n\}. \ \ \ \ \ \ \ \ \ \ (3.1)$$
 In  general, the analysis of  solutions of the system (3.1) is difficult. Therefore we shall consider some particular solutions of the system (3.1):
 \\ 1) Trivial solution $(0,...,0) .$
 \\ 2) Assume that $a_{k\pi(k)} + a_{j_{k}\pi(k)}=d$ for any $k=1,...,n$ , where $d=const,\ d\neq 0,$ then one of the particular solutions
 is $x_{i}=\frac{1}{d}, \ i=1,...,n.$

 Now we consider the system of equations for $n=2$, then permutations (with $\pi\neq\tau$) will be
 $$ \pi=\left( \begin{array}{cc} 1  \ \  \ 2 \ \\
 2 \ \ \  1 \ \end{array} \right),\ \ \
    \tau=\left( \begin{array}{cc} 1  \ \  \ 2 \ \\
 1 \ \ \  2 \ \end{array} \right).$$
 Then system of equations (3.1) will be

 $$\left\{\begin{array}{ll} a_{12}x^{2}_{1}+a_{22}x^{2}_{2}=x_{2} \\

   a_{21}x^{2}_{2}+a_{11}x^{2}_{1}=x_{1}.
   \end{array} \right.$$
 Assume that \ $a_{ij}\neq0 \ ; \ i,j\in\{1,2\}$. For easiness let's denote $a_{12}=a ,\  a_{22}=b ,\ a_{21}=c,\  a_{11}=d $ and $ x_{1}=x,\ x_{2}=y $ then we have \ \ \ \ $$\left\{\begin{array}{ll} ax^{2} + by^{2}=y \\
   dx^{2} + cy^{2}=x
   \end{array} \right.\ \  $$
 from this system we have
 $$ (bd-ac)^{2}x^{4} - 2b(bd-ac)x^{3}+ (b^{2} + cd)x^{2} - cx = 0 \ \ \ \ \ \ (3.2)$$
 By full analysis of (3.2) we get

 \

 \textbf{Theorem 3.1.} \emph{Solutions of the equation (3.2) are $x=0$ and
\\ 1) If $bd=ac$, $b^{2}+cd\neq 0$ then $ x=\frac{c}{b^{2}+cd}$  .
\\ 2) Let $p=(3cd - b^{2})/3(bd - ac)^{2}, \ q=2(9bcd + b^{3})/27(bd-ac)^{3} - c/(bd - ac)^{2},$
\\  $\Delta=(q/2)^{2}+ (p/3)^{3}$ if $bd\neq ac$, then
\\ a) for $\Delta <0$, then there exist three real solutions.
\\ b) for $\Delta >0$, then there is one  real solution and two complex conjugate solutions.
\\ c) for $\Delta =0$, and $p\neq0, q\neq0$ then there are two  real solutions.
\\ d) for $\Delta =0$, and $p=q=0 $ then there exist one real solution i.e. $x= 2b/3(bd-ac)$.}\\

\textbf{4. Isomorphism of $ E^{n}_{\pi,\tau }$.}

\

 \textbf{Proposition 4.1.}\ \emph{Let $ E^{n}_{\pi,\tau }$ be an evolution algebra of permutations with the following conditions :\\
 1) $a_{i\pi(i)}\cdot a_{i\tau(i)}\neq 0, \ 1\leq i \leq n,$  }\\
 2) $\pi=\pi_{1}\circ \pi_{2}\circ ... \circ \pi_{r}, \ \tau=\tau_{1}\circ \tau_{2}\circ ... \circ\tau_{r}, \ where \
 \pi_{1}=(\pi_{11}, \pi_{12}, ... ,\pi_{1k_{1}}), \ \pi_{2}=(\pi_{21}, \pi_{22}, ... ,\pi_{2k_{2}}),\ ... ,\ \pi_{r}=(\pi_{r1}, \pi_{r2}, ... ,\pi_{rk_{r}}) \ and \ \  \tau_{1}=(\tau_{11},\tau_{12}, ... ,\tau_{1k_{1}}), \tau_{2}=(\tau_{21},\tau_{22},...,\tau_{2k_{2}}),\\ ...,\ \tau_{r}=(\tau_{r1},\tau_{r2},...,\tau_{rk_{r}}) \ are \ independent\  cycles\  of \  \pi \ and \  \tau \ $
 $ respectively,\  and \ k_{1}+k_{2}+...+k_{r}=n, \ \tau_{im}\in \{\pi_{i1},\pi_{i2},...,\pi_{ik_{i}}\},\ 1\leq i \leq r,\ 1\leq m \leq k_{i}\
 (i.e.\  \pi_{k}\  and \ \tau_{k}\ consist \ of \ one\ and\\  the\ same \ elements,\ only\ has\
  difference\ between\ seats\ of\ elements ).\ $ \\
 $Then$
\ \ \ \ \ \ \ \ \ \ $$E^{n}_{\pi,\tau }\cong E^{k_{1}}_{\pi_{1},\tau_{1}}\oplus E^{k_{2}}_{\pi_{2},\tau_{2}}\oplus ...\oplus E^{k_{r}}_{\pi_{r},\tau_{r}}. $$

\

\textbf{Proof}. $E^{k_{i}}_{\pi_{i},\tau_{i}}$ is $k_{i}$-dimensional evolution algebra of permutations $\pi_{i}$ and $\tau_{i}$. Indeed, $$\left\{\begin{array}{ll} e_{\pi_{im}}\cdot e_{\pi_{im}}=a_{\pi_{im}\pi(\pi_{im})}\cdot e_{\pi(\pi_{im})} + a_{\pi_{im}\tau(\pi_{im})}\cdot e_{\tau(\pi_{im})} \\
e_{\pi_{im}}\cdot e_{\pi_{ik}}=0,\ \ m\neq k
   \end{array} \right.\ \  $$
and $\pi_{im},\pi(\pi_{im})=\pi_{i,m+1}\in \{\pi_{i1},\pi_{i2},...,\pi_{ik_{i}}\},\ \tau(\pi_{im})=\tau(\tau_{is})\in \{\tau_{i1},\tau_{i2},...,\tau_{ik_{i}}\},\\  1\leq m,s \leq k_{i}.$

The isomorphism is provided by the following change of basis $$e_{im}=e_{\pi_{im}}( \ or \ e_{im}=e_{\tau_{im}} ),\ 1\leq i \leq r,\ 1\leq m \leq k_{i}. $$
Thus we have the evolution algebra $ E^{k_{i}}_{\pi_{i},\tau_{i}}$ with the basis $e_{im},\ 1\leq i \leq r,\ 1\leq m \leq k_{i} $ .

If $i\neq j$ then $\pi_{im}\neq \pi_{jk}$ or $\tau_{im}\neq \tau_{jk}$ for any $m$ and $k$, it means
$$e_{\pi_{im}}\neq e_{\pi_{jk}},\  e_{\pi_{im}}\cdot e_{\pi_{jk}}=0.$$
It follows that $E^{k_{i}}_{\pi_{i},\tau_{i}} \cap E^{k_{j}}_{\pi_{j},\tau_{j}}=\emptyset$ (i.e $\{ \pi_{i1},\pi_{i2},...,\pi_{ik_{i}} \} \cap \{ \pi_{j1},\pi_{j2},...,\pi_{jk_{j}} \}=\emptyset $ ) for $i\neq j$
and we have

$$E^{n}_{\pi,\tau }\cong E^{k_{1}}_{\pi_{1},\tau_{1}}\oplus E^{k_{2}}_{\pi_{2},\tau_{2}}\oplus ...\oplus E^{k_{r}}_{\pi_{r},\tau_{r}}. $$

\

\textbf{Proposition 4.2.}\ \emph{Any evolution algebra of permutations $E^{n}_{\pi,\tau_{0} }$ with permutations
$$\pi=(k_{1},k_{2},...,k_{n}),\ \tau_{0}=(1)(2) ... (n)\  and \  a_{i\pi(i)}\cdot a_{i\tau(i)}\neq 0, \ 1\leq i \leq n $$
is isomorphic to the evolution algebra  $E^{n}_{\pi_{\ast},\tau_{0}} $, with the table of multiplications given by:
$$ \left\{\begin{array}{ll} e^{'}_{i}\cdot e^{'}_{i}= a_{\pi ^{i-1}(1)\pi^{i}(1)}e^{'}_{i+1} + a_{\pi ^{i-1}(1)\pi^{i-1}(1)}e^{'}_{i}, \ \ 1\leq i \leq n \\
   e^{'}_{i}\cdot e^{'}_{j}=0,\ \ i\neq j
   \end{array} \right.\ \  $$
   where  $\pi_{\ast}=(1,2,...,n), \ $ and $a_{\pi ^{i-1}(1)\pi^{i}(1)}\cdot a_{\pi ^{i-1}(1)\pi^{i-1}(1)}\neq 0
   , \ \pi^{0}(1)=1 $   }

\

\textbf{Proof}. The isomorphism is established by basis permutations:
$$ e^{'}_{1}=e_{1}, \ e^{'}_{i}=e_{\pi^{i-1}(1)}, \ \ 2\leq i \leq n. $$
Indeed, according to $\tau_{1}(i)=i $ we get

$ e^{'}_{i}\cdot e^{'}_{i}=e_{\pi^{i-1}(1)}\cdot e_{\pi^{i-1}(1)}=a_{\pi ^{i-1}(1)\pi^{i}(1)} \cdot e_{\pi^{i}(1)} + a_{\pi ^{i-1}(1)\tau ( \pi^{i-1}(1))}\cdot e_{\tau ( \pi^{i-1}(1))}= a_{\pi ^{i-1}(1)\pi^{i}(1)} \cdot e_{\pi^{i}(1)} + a_{\pi ^{i-1}(1) \pi^{i-1}(1)}\cdot e_{ \pi^{i-1}(1)}=a_{\pi ^{i-1}(1)\pi^{i}(1)}e^{'}_{i+1} + a_{\pi ^{i-1}(1)\pi^{i-1}(1)}e^{'}_{i}.$

\

\textbf{Proposition 4.3.}\ \emph{ Any evolution algebra of permutations $E^{n}_{\pi,\tau }$ with permutations
$\pi=(k_{1},k_{2},...,k_{n}),\ \tau=(l_{1},l_{2},...,l_{n})$ and $\pi\cdot \tau=(1)(2) ... (n)$ is isomorphic to the evolution algebra $E^{n}_{\pi_{\ast},\tau_{\ast}} $,
 which a table of multiplications given by:
  $$ \left\{\begin{array}{ll} e^{'}_{i}\cdot e^{'}_{i}= a_{\pi ^{i-1}(1)\pi^{i}(1)}e^{'}_{i+1} + a_{\pi ^{i-1}(1)\pi^{i-2}(1)}e^{'}_{i-1}, \ \ 1\leq i \leq n \\
   e^{'}_{i}\cdot e^{'}_{j}=0,\ \ i\neq j
   \end{array} \right.\ \  $$
 with permutations $\pi_{\ast}=(1,2,...,n), \ \tau_{\ast}=(1,n,n-1,...,2)$ and $a_{\pi ^{i-1}(1)\pi^{i}(1)}\cdot a_{\pi ^{i-1}(1)\pi^{i-2}(1)}\neq 0. $}

 \

\textbf{Proof}. Similarly to the proof of Proposition 4.2 taking the change of basis
$$ e^{'}_{1}=e_{1}, \ e^{'}_{i}=e_{\pi^{i-1}(1)}, \ \ 2\leq i \leq n. $$ and we have $\pi\cdot \tau=(1)(2) ... (n)$ it means that
$\tau(\pi (k))=k,\ \tau(\pi^{i+1}(1))=\tau(\pi(\pi^{i}(1)))=\pi^{i}(1)$ thus we have

$ e^{'}_{i}\cdot e^{'}_{i}=e_{\pi^{i-1}(1)}\cdot e_{\pi^{i-1}(1)}=a_{\pi ^{i-1}(1)\pi^{i}(1)} \cdot e_{\pi^{i}(1)} + a_{\pi ^{i-1}(1)\tau ( \pi^{i-1}(1))}\cdot e_{\tau ( \pi^{i-1}(1))}= a_{\pi ^{i-1}(1)\pi^{i}(1)} \cdot e_{\pi^{i}(1)} + a_{\pi ^{i-1}(1) \pi^{i-2}(1)}\cdot e_{ \pi^{i-2}(1)}=a_{\pi ^{i-1}(1)\pi^{i}(1)}e^{'}_{i+1} + a_{\pi ^{i-1}(1)\pi^{i-2}(1)}e^{'}_{i-1}. $

\

\end{document}